\documentclass[12pt,reqno,letterpaper]{amsart}

\DeclareMathOperator{\ann}{ann}%
\DeclareMathOperator{\chr}{char}%
\DeclareMathOperator{\proj}{proj}%

\usepackage{amscd,amsmath,amssymb,latexsym}
\usepackage{hyperref}

\begin{document}

\newtheorem*{theorem*}{Theorem}
\newtheorem{theorem}{Theorem}
\newtheorem{corollary}[theorem]{Corollary}
\newtheorem{proposition}{Proposition}

\theoremstyle{definition}
\newtheorem*{definition*}{Definition}
\newtheorem*{remark*}{Remark}

\newcommand{\Ac}{\mathcal{A}}
\newcommand{\Fix}{\mathrm{Fix}}
\newcommand{\Fp}{\mathbb{F}_p}
\newcommand{\Gal}{\mathrm{Gal}}
\newcommand{\N}{\mathbb{N}}
\newcommand{\Z}{\mathbb{Z}}

\parskip=10pt minus 2pt

\title[Automatic Realizations of Galois Groups]
{Automatic Realizations of Galois Groups with Cyclic Quotient of
Order $p^n$}

\author[J\'{a}n Min\'{a}\v{c}]{J\'an Min\'a\v{c}$^\dag$}
\address{Department of Mathematics, Middlesex College, \ University
of Western Ontario, London, Ontario \ N6A 5B7 \ CANADA}
\thanks{$^\dag$Research supported in part by NSERC grant R0370A01,
and by a Distinguished Research Professorship at the University of
Western Ontario.}
\email{minac@uwo.ca}

\author[Andrew Schultz]{Andrew Schultz}
\address{Department of Mathematics, Building 380, Stanford
University, Stanford, California \ 94305-2125 \ USA}
\email{aschultz@stanford.edu}

\author[John Swallow]{John Swallow$^\ddag$}
\address{Department of Mathematics, Davidson College, Box 7046,
Davidson, North Carolina \ 28035-7046 \ USA}
\thanks{$^\ddag$Research supported in part by National Security
Agency grant MDA904-02-1-0061.}
\email{joswallow@davidson.edu}


\begin{abstract}
    We establish automatic realizations of Galois groups among
    groups $M\rtimes G$, where $G$ is a cyclic group of order $p^n$
    for a prime $p$ and $M$ is a quotient of the group ring
    $\Fp[G]$.
\end{abstract}

\date{March 24, 2006}

\maketitle

The fundamental problem in inverse Galois theory is to determine,
for a given field $F$ and a given profinite group $G$, whether there
exists a Galois extension $K/F$ such that $\Gal(K/F)$ is isomorphic
to $G$. A natural sort of reduction theorem for this problem takes
the form of a pair $(A,B)$ of profinite groups with the property
that, for all fields $F$, the existence of $A$ as a Galois group
over $F$ implies the existence of $B$ as a Galois group over $F$. We
call such a pair an automatic realization of Galois groups and
denote it $A\implies B$.  The trivial automatic realizations are
those given by quotients of Galois groups; by Galois theory, if $G$
is realizable over $F$ then so is every quotient $H$.  It is a
nontrivial fact, however, that there exist nontrivial automatic
realizations.  (See \cite{Je1,Je2,Je3} for a good overview of the
theory of automatic realizations.  Some interesting automatic
realizations of groups of order $16$ are obtained in \cite{GrSm},
and these and other automatic realizations of finite 2-groups are
collected in \cite{GSS}.  For comprehensive treatments of related
Galois embedding problems, see \cite{JLY} and \cite{Le}.)

The usual techniques for obtaining automatic realizations of Galois
groups involve an analysis of Galois embedding problems. In this
paper we offer a new approach based on the structure of natural
Galois modules: we use equivariant Kummer theory to reformulate
realization problems in terms of Galois modules, and then we solve
Galois module problems.  We take this approach in proving
Theorem~\ref{th:main}, which establishes automatic realizations for
a useful family of finite metacyclic $p$-groups.  Our methods extend
those of \cite{MS2} and \cite{MSS}.  It is interesting to observe
that, although not visible here, the essential fact underpinning our
results is Hilbert~90.  Indeed, the structural results in \cite{MSS}
rely crucially on the repeated application of Hilbert 90, using
combinatorial and Galois-theoretic arguments to draw out the
consequences.

Let $p$ be a prime, $n\in \N$, and $G$ a cyclic group of order $p^n$
with generator $\sigma$.  For the group ring $\Fp[G]$, there exist
precisely $p^n$ nonzero ring quotients, namely $M_j :=
\Fp[G]/\langle (\sigma-1)^j\rangle$ for $j=1, 2, \dots, p^n$.
Multiplication in $\Fp[G]$ induces an $\Fp[G]$-action on each $M_j$.
In particular, each $M_j$ is a $G$-module. Let $M_j\rtimes G$ denote
the semidirect product.
\begin{theorem}\label{th:main}
    We have the following automatic realizations of Galois groups:
    \begin{equation*}
        M_{p^i+c}\rtimes G \implies M_{p^{i+1}}\rtimes G, \quad\ \
        0\le i<n, \ \ 1\le c<p^{i+1}-p^i.
    \end{equation*}
\end{theorem}

In Section~\ref{se:gm} we recall some facts about the set of
quotients $M_j$ and the semidirect products $M_j\rtimes G$.  In
Sections~\ref{se:kummer} through \ref{se:jepsilon} we consider the
case $\chr F\neq p$.  Following Waterhouse \cite{W}, we recall in
Section~\ref{se:kummer} a generalized Kummer correspondence over
$K$, where $K$ is a cyclic extension of $F$ of degree $p^n$, and in
Section~\ref{se:index} we establish a proposition detecting when
such extensions are Galois over $F$. In Section~\ref{se:jepsilon} we
decompose $J_\epsilon$, the crucial Kummer submodule of the module
$K(\xi_p)^\times / K(\xi_p)^{\times p}$, as an $\Fp[\Gal(K(\xi_p)/
F(\xi_p))]$-module, where $\xi_p$ is a primitive $p$th root of
unity. In Section~\ref{se:proof} we prove Theorem~\ref{th:main},
using Sections~\ref{se:kummer} through \ref{se:jepsilon} in the case
$\chr F\neq p$ and Witt's Theorem in the case $\chr F=p$. The case
$i=0$ was previously considered by two of the authors \cite[Theorem
1(A)]{MS2}.

\section{Groups and $\Fp[G]$-modules}\label{se:gm}

Let $p$ be a prime and $G=\langle \sigma\rangle$ an abstract group
of order $p^n$.  We recall some facts concerning $R$-modules, where
$R$ is the group ring $\Fp[G]$.  Because we frequently view $R$ as a
module over $R$, to prevent confusion we write the module $R$ as
\begin{equation*}
    R=\oplus_{j=0}^{p^n-1} \Fp \tau^j,
\end{equation*}
where $\sigma$ acts by multiplication by $\tau$. For convenience we
set $\rho:=\sigma-1$.

The set of nonzero cyclic $R$-modules is identical to the set of
nonzero indecomposable $R$-modules, and these are precisely the
$p^n$ quotients $M_j:=R/\langle (\tau-1)^j\rangle$, $1\le j\le p^n$.
Each $M_j$ is a local ring, with unique maximal ideal $\rho M_j$,
and is annihilated by $\rho^{j}$ but not $\rho^{j-1}$.

Moreover, for each $j$ there exists a $G$-equivariant isomorphism
from $M_j$ to its dual $M_j^*$, as follows.  For each $i\in \{1,
\dots, p^n\}$ we choose the $\Fp$-basis of $M_j$ consisting of the
images of $\{1,(\tau-1),\dots,(\tau-1)^{j-1}\}$ and define an
$\Fp$-linear map $\lambda:M_j\to \Fp$ by
\begin{equation*}
    \lambda\left(f_0+f_1\overline{(\tau-1)}+\cdots+f_{j-1}
    \overline{(\tau-1)}^{j-1}\right) = f_{j-1},
\end{equation*}
where $f_k \in \Fp$, $k=0, \dots, j-1$.  Observe that $\ker\lambda$
contains no nonzero ideal of $M_j$.  Then
\begin{equation*}
    Q:M_j\times M_j\to \Fp, \qquad Q(a,b):=\lambda(ab),\ a,b\in M_j
\end{equation*}
is a nonsingular symmetric bilinear form.  Thus $M_j$ is a symmetric
algebra. (See \cite[page~442]{La}.)  Moreover, $Q$ induces a
$G$-equivariant isomorphism $\psi:M_j\to M_j^*$ given by
$(\psi(a))(b) = Q(a,b)$, $a,b \in M_j$.

\begin{remark*}
    In order for $\psi$ to be $G$-equivariant, we must define the
    action on $M_j^*$ by $\sigma f(m) = f(\sigma m)$ for all $m \in
    M_j$, and since $G$ is commutative, this action is well-defined.
    It is worthwhile to observe, however, that $M_j^*$ is
    $\Fp[G]$-isomorphic to the module $\tilde M_j^*$ on which the
    action of $G$ is defined by $\sigma f(m) = f(\sigma^{-1} m)$ for
    all $m \in M_j$. Indeed by the $G$-equivariant isomorphism between
    $M_j$ and $M_j^*$ it is sufficient to show that the
    $\Fp[G]$-module $\tilde M_j$ obtained from $M_j$ by twisting the
    action of $G$ via the automorphism $\sigma \to \sigma^{-1}$ is
    naturally isomorphic to $M_j$. But this follows readily by
    extending the automorphism $\sigma \to \sigma^{-1}$ to the
    automorphism of the group ring $\Fp[G]$ and then inducing the
    required $\Fp[G]$-isomorphism between $M_j$ and $M_j^*$.
\end{remark*}

We also recall some facts about the semidirect products $H_j :=
M_j\rtimes G$, $j=1,\dots, p^n$.  For each $j$, the group $H_j$ has
order $p^{j+n}$; exponent $p^n$, except when $j=p^n$, in which case
the exponent is $p^{n+1}$; nilpotent index $j$; rank (the smallest
number of generators) $2$; and Frattini subgroup $\Phi(H_j) =(\rho
M_j)\rtimes G^p$.  Finally, for $j<k$, $H_j$ is a quotient of $H_k$
by the normal subgroup $\rho^{j}M_k\rtimes 1$.

\section{Kummer Theory with Operators}\label{se:kummer}

For Sections~\ref{se:kummer} through \ref{se:jepsilon} we adopt the
following hypotheses.  Suppose that $G=\Gal(K/F)=\langle
\sigma\rangle$ for an extension $K/F$ of degree $p^n$ of fields of
characteristic not $p$. We let $\xi_p$ be a primitive $p$th root of
unity and set $\hat F:=F(\xi_p)$, $\hat K:=K(\xi_p)$, and $J:=\hat
K^\times/\hat K^{\times p}$, where $\hat K^\times$ denotes the
multiplicative group $\hat K\setminus \{0\}$. We write the elements
of $J$ as $[\gamma]$, $\gamma\in \hat K^\times$, and we write the
elements of $\hat F^\times/\hat F^{\times p}$ as $[\gamma]_{\hat
F}$, $\gamma\in \hat F^\times$. We moreover let $\epsilon$ denote a
generator of $\Gal(\hat F/F)$ and set $s=[\hat F:F]$. Since $p$ and
$s$ are relatively prime, $\Gal(\hat K/F)\simeq \Gal(\hat F/F)\times
\Gal(K/F)$. Therefore we may naturally extend $\epsilon$ and
$\sigma$ to $\hat K$, and the two automorphisms commute in
$\Gal(\hat K/F)$.  Using the extension of $\sigma$ to $\hat K$, we
write $G$ for $\Gal(\hat K/\hat F)$ as well.  Then $J$ is an
$\Fp[G]$-module.  Finally, we let $t\in \Z$ such that
$\epsilon(\xi_p)=\xi_p^t$.  Then $t$ is relatively prime to $p$, and
we let $J_\epsilon$ be the $t$-eigenspace of $J$ under the action of
$\epsilon$: $J_\epsilon = \{ [\gamma] \ :\  \epsilon [\gamma] =
[\gamma]^t \}$.

Observe that since $\epsilon$ and $\sigma$ commute, $J_\epsilon$ is
an $\Fp[G]$-subspace of $J$.  By \cite[\S 5, Proposition]{W}, we
have a Kummer correspondence over $K$ of finite subspaces $M$ of the
$\Fp$-vector space $J_\epsilon$ and finite abelian exponent $p$
extensions $L$ of $K$:
\begin{multline*}
    M= ((\hat K L)^{\times p}\cap \hat K^\times)/\hat K^{\times p}
    \leftrightarrow \\ L=L_M = \text{maximal $p$-extension of $K$
    in\ } \hat L_{M} := \hat K(\root{p}\of{\gamma}:[\gamma]\in M).
\end{multline*}
As Waterhouse shows, for $M\subset J_\epsilon$, the automorphism
$\epsilon\in \Gal(\hat K/K)$ has a unique lift $\tilde\epsilon$ to
$\Gal(\hat L_M/K)$ of order $s$, and $L_M$ is the fixed field of
$\tilde\epsilon$.

In the next proposition we provide some information about the
corresponding Galois modules when $L_M/F$ is Galois.  Recall that in
the situation above, the Galois groups $\Gal (L_M/F)$ and $\Gal(\hat
L_M/\hat K)$ are naturally $G$-modules under the action induced by
conjugations of lifts of the elements in $G$ to $\Gal(L_M/F)$ and
$\Gal (\hat L_M/\hat F)$. Furthermore, because the Galois groups
have exponents dividing $p$, we see that $\Gal(L_M/F)$ and
$\Gal(\hat L_M/\hat F)$ are in fact $\Fp[G]$-modules.

\begin{proposition}\label{pr:kummer}
    Suppose that $M$ is a finite $\Fp$-subspace of $J_\epsilon$.  Then
    \begin{enumerate}
        \item\label{it:k1} $L_M$ is Galois over $F$ if and only if
        $M$ is an $\Fp[G]$-submodule of $J_\epsilon$.
        \item\label{it:k2} If $L_M/F$ is Galois, then base extension
        $F\to \hat F$ induces a natural isomorphism of $G$-modules
        $\Gal(L_M/F)\simeq \Gal(\hat L_M/\hat F)$.
        \item\label{it:k3} If $L_M/F$ is Galois, then as
        $G$-modules,
        \begin{equation*}
            \Gal(L_M/K)\simeq \Gal(\hat L_M/\hat K)\simeq M.
        \end{equation*}
    \end{enumerate}
\end{proposition}

\begin{proof}
    (\ref{it:k1}). Suppose first that $L_M/F$ is Galois.  Then $\hat
    L_{M}=L\hat K/\hat F$ is Galois as well.  Every automorphism of
    $\hat K$ extends to an automorphism of $\hat L_M$, and therefore
    $M$ is an $\Fp[G]$-submodule of $J$.  From \cite[\S 5,
    Proposition]{W} we see that $M$ is an $\Fp[G]$-submodule of
    $J_\epsilon$.

    Going the other way, suppose that $M$ is a finite
    $\Fp[G]$-submodule of $J_\epsilon$.  By the correspondence
    above, $L_M/K$ is Galois. Then $M$ is also an $\Fp[\Gal(\hat
    K/F)]$-submodule of $J_\epsilon$ and therefore $\hat L_M/F$ is
    Galois.  Now since $K/F$ is Galois, every automorphism of $\hat
    L_M$ sends $K$ to $K$.  Moreover, since $L_M$ is the unique
    maximal $p$-extension of $K$ in $\hat L_M$, every automorphism
    of $\hat L_{M}$ sends $L_M$ to $L_M$.  Therefore $L_M/F$ is
    Galois.

    (\ref{it:k2}). Suppose $L_M/F$ is Galois. Since $\hat F/F$ and
    $L_M/F$ are of relatively prime degrees, we have $\Gal(L_M\hat
    F/F)\simeq \Gal(\hat F/F)\times \Gal(L_M/F)$.  Therefore we have
    a natural isomorphism $G=\Gal(K/F)\simeq \Gal(\hat K/\hat F)$,
    and the natural isomorphism $\Gal(\hat L_M/\hat
    F)\simeq\Gal(L_M/F)$ is an isomorphism of $G$-extensions.

    (\ref{it:k3}).  Suppose $L_M/F$ is Galois.  By (\ref{it:k2}), it
    is enough to show that $\Gal(\hat L_M/\hat K)\simeq M$ as
    $G$-modules. Under the standard Kummer correspondence over $\hat
    K$, finite subspaces of the $\Fp$-vector space $J$ correspond to
    finite abelian exponent $p$ extensions $\hat L_M$ of $\hat K$,
    and $M$ and $\Gal(\hat L_M/\hat K)$ are dual $G$-modules under a
    $G$-equivariant canonical duality $\langle m,g\rangle =
    g(\root{p}\of{m})/\root{p}\of{m}$.  (See \cite[pages 134 and
    135]{W} and \cite[\S 2.3]{MS2}.) Because $M$ is finite, $M$
    decomposes into a direct sum of indecomposable $\Fp[G]$-modules.
    From Section~\ref{se:gm}, all indecomposable $\Fp[G]$-modules
    are $G$-equivariant self-dual modules.  Hence there is a
    $G$-equivariant isomorphism between $M$ and its dual $M^*$, and
    $\Gal(\hat L_M/\hat F)\simeq M$ as $G$-modules.
\end{proof}

\section{The Index}\label{se:index}

We keep the same assumptions given at the beginning of
Section~\ref{se:kummer}.  Set $A := \ann_J \rho^{p^n-1} = \{
[\gamma] \in J  :  \rho^{p^{n-1}-1} [\gamma] = [1] \}$. The
following homomorphism appears in a somewhat different form in
\cite[Theorem 3]{W}:
\begin{definition*}
    The \emph{index} $e(\gamma)\in \Fp$ for
    $[\gamma]\in A$ is defined by
    \begin{equation*}
        \xi_p^{e(\gamma)} = \left( \root{p} \of{N_{\hat K/\hat
        F}(\gamma)}\right)^{\rho}.
    \end{equation*}
\end{definition*}
The index is well-defined, as follows.  First, since
\begin{equation*}
    1+\sigma+\dots+\sigma^{p^n-1}=(\sigma-1)^{p^n-1}=\rho^{p^n-1}
\end{equation*}
in $\Fp[G]$, $[N_{\hat K/\hat F}(\gamma)]=[\gamma]^{\rho^{p^n-1}}$,
which is the trivial class $[1]$ by the assumption $[\gamma]\in A$.
As a result, $\root{p}\of{N_{\hat K/\hat F}(\gamma)}$ lies in $\hat
K$ and is acted upon by $\sigma$ and therefore $\rho$. Observe that
$e(\gamma)$ depends neither on the representative $\gamma$ of
$[\gamma]$ nor on the particular $p$th root of $N_{\hat K/\hat
F}(\gamma)$.

The index function $e$ is a group homomorphism from $A$ to $\Fp$.
Therefore the restriction of $e$ to any submodule of $A$ is either
trivial or surjective.  Moreover, the index is trivial for any
$[\gamma]$ in the image of $\rho$:
\begin{equation*}
    \xi_p^{e(\gamma^\rho)} = \root{p}\of{ N_{\hat K/\hat
    F}(\gamma^\rho)}^{\rho} = \root{p}\of{(N_{\hat K/\hat
    F}(\gamma))^\rho} = \root{p}\of{1}^{\rho} = 1,
\end{equation*}
or $e(\gamma^\rho)=0$.

Following Waterhouse, we show how the index function permits the
determination of $\Gal(\hat L_M/\hat F)$ as a $G$-extension.

For $1\le j\le p^n$ and $e\in \Fp$, write $H_{j,e}$ for the group
extension of $M_j$ by $G$ with $\tilde \sigma^{p^n} =
e(\tau-1)^{j-1}$, where $\tilde\sigma$ is a lift of $\sigma$.
Observe that $H_{j,0}=H_j=M_j\rtimes G$.

Let $N_\gamma$ denote the cyclic $\Fp[G]$-submodule of $J$ generated
by $[\gamma]$.
\begin{proposition}\label{pr:groupoflm} (See \cite[Theorem 2]{W}.)
    Let $[\gamma]\in J_\epsilon$ and $M=N_\gamma$.
    \begin{enumerate}
        \item\label{it:lm1} If $M\simeq M_j$ for $1\le j<p^n$ and
        $e=e(\gamma)$, then $\Gal(L_M/F)\simeq H_{j,e}$ as
        $G$-extensions.
        \item\label{it:lm2} If $M\simeq \Fp[G]$ then
        $\Gal(L_M/F)\simeq \Fp[G]\rtimes G$.
    \end{enumerate}
\end{proposition}
Before presenting the proof, we note that Waterhouse also tells us
that for $j<p^n$ and $e\neq 0$, $H_{j,e}\not\simeq H_j$ (see
\cite[Theorem 2]{W}), and for $j=p^n$ then there is a $G$-extension
isomorphism $H_{p^n,e}\simeq H_{p^n}$ for every $e$. In particular,
we may use Proposition~\ref{pr:groupoflm} later to deduce that if
$M\simeq M_j$ for $j<p^n$ and $\Gal(L_M/F)\simeq M_j\rtimes G$, then
$e(\gamma)=0$.
\begin{proof}
    Suppose $M\simeq M_j$ for some $1\le j\le p^n$. By
    Proposition~\ref{pr:kummer}(\ref{it:k3}), $\Gal(L_M/K)\simeq
    M_j$ as $G$-modules.  Hence $\Gal(L_M/F)\simeq H_{j,e}$ for some
    $e$. If $j=p^n$ then from the isomorphism $H_{p^n,e}\simeq
    H_{p^n}$ above we have the second item. By
    Proposition~\ref{pr:kummer}(\ref{it:k2}), it remains only to
    show that if $j<p^n$, $\Gal(\hat L_M/\hat F)\simeq
    H_{j,e(\gamma)}$.

    Let $\tilde\sigma$ denote a pullback of $\sigma\in G$ to
    $\Gal(\hat L_M/\hat F)$. Then $\tilde\sigma^{p^n}$ lies in
    $Z(\Gal(\hat L_M/\hat F)) \cap \Gal(\hat L_M/\hat K)$, where
    $Z(\Gal(\hat L_M/\hat F))$ denotes the center of $\Gal(\hat
    L_M/\hat F).$  Using the $G$-equivariant Kummer pairing
    \begin{equation*}
        \langle \cdot,\cdot\rangle \colon \Gal(\hat L_M/\hat K)
        \times M \to \langle \xi_p \rangle \simeq \Fp
    \end{equation*}
    we see that $Z(\Gal (\hat L_M/\hat K))$ annihilates $\rho M$.
    Furthermore, since this pairing is nonsingular we deduce that
    $Z(\Gal(\hat L_M/\hat K)) \simeq M/\rho M$ and we can choose a
    generator $\eta$ of $Z(\Gal (\hat L_M/\hat K))$ such that
    \begin{equation*}
        \langle \eta, [\gamma] \rangle =
        \eta(\root{p}\of{\gamma})/\root{p}\of{\gamma} = \xi_p.
    \end{equation*}
    In particular, if $\tilde \sigma^{p^n} = \eta^e$ then
    \begin{equation*}
        (\root{p}\of{\gamma})^{(\tilde\sigma^{p^n}-1)} =
        \xi_{p}^{e}.
    \end{equation*}
    Therefore
    \begin{equation*}
        \root{p}\of{\gamma}^{(\tilde\sigma^{p^n}-1)} =
        \root{p}\of{\gamma}^{(1+\tilde\sigma+
        \dots+\tilde\sigma^{p^n-1}) (\tilde\sigma-1)} =
        \left(\root{p}\of{ N_{\hat K/\hat F}(\gamma)}\right)^{\rho}=
        \xi_p^{e(\gamma)}.
    \end{equation*}
\end{proof}

\section{The $\Fp[G]$-module $J_\epsilon$}\label{se:jepsilon}

Again we keep the same assumptions given at the beginning of
Section~\ref{se:kummer}.  In this section we develop the crucial
technical results needed for Theorem~\ref{th:main}: a decomposition
of the $\Fp[G]$-module $J_\epsilon$ into cyclic direct summands, and
a determination of the value of the index function $e$ on certain of
the summands.

We first show that $J_\epsilon$ is indeed a summand of $J$.  Then we
combine a decomposition of $J$ into indecomposables, taken from
\cite[Theorem 2]{MSS}, with uniqueness of decompositions into
indecomposables, to achieve important restrictions on the possible
summands of $J_\epsilon$.  Much of the remainder of the proof is
devoted to establishing that we have an ``exceptional summand'' of
dimension $p^r+1$ on which the index function is nontrivial.  In the
argument we need \cite[Proposition 7]{MSS} in particular to derive a
lower bound for the dimension of that summand.

\begin{theorem}\label{th:jepsilon}
    Suppose that $p>2$ or $n>1$.  The $\Fp[G]$-module $J_\epsilon$
    decomposes into a direct sum $J_\epsilon = U\oplus_{a\in \Ac}
    V_\alpha$, with $\Ac$ possibly empty, with the following
    properties:
    \begin{enumerate}
        \item\label{it:j1} For each $\alpha\in \Ac$ there exists
        $i\in \{0,\dots,n\}$ such that $V_\alpha\simeq M_{p^i}$.
        \item\label{it:j2} $U\simeq M_{p^r+1}$ for some $r\in
        \{-\infty,0,1,\dots,n-1\}$.
        \item\label{it:j3} $e(U)=\Fp$.
        \item\label{it:j4} If $V_\alpha\simeq M_{p^i}$ for $0\le
        i\le r$, then $e(V_\alpha)=\{0\}$.
    \end{enumerate}
\end{theorem}
\noindent Here we observe the convention that $p^{-\infty}=0$.
\begin{proof}
    We show first that $J_\epsilon$ is a direct summand of $J$ by
    adapting an approach to descent from \cite[page~258]{Sa}.  Recall
    that $[\hat F:F] = s$ and $\epsilon(\xi_p)= \xi_p^t$. Thus $s$
    and $t$ are both relatively prime to $p$. Let $z\in \Z$ satisfy
    $zst^{s-1}\equiv 1\ (\bmod\ p)$, and set
    \begin{equation*}
        T = z\cdot \sum_{i=1}^st^{s-i}\epsilon^{i-1} \in
        \Z[\Gal(\hat K/F)].
    \end{equation*}
    We calculate that $(t-\epsilon)T\equiv 0\ (\bmod\ p)$, and hence
    the image of $T$ on $J$ lies in $J_\epsilon$.  Moreover,
    $\epsilon$ acts on $J_\epsilon$ by multiplication by $t$, and
    therefore $T$ acts as the identity on $J_\epsilon$.  Finally,
    since $\epsilon$ and $\sigma$ commute, $T$ and $I-T$ commute
    with $\sigma$. Hence $J$ decomposes into a direct sum
    $J_\epsilon\oplus J_\nu$, with associated projections $T$ and
    $I-T$.

    We claim that $e((I-T)A)=\{0\}$.  Since $\xi_p\in \hat F$, the
    fixed field $\Fix_{\hat K}(\sigma^p)$ may be written $\hat
    F(\root{p}\of{a})$ for a suitable $a\in \hat F^\times$.  By
    \cite[\S 5, Proposition]{W}, $\epsilon([a]_{\hat F})=[a]_{\hat
    F}^t$.  Suppose $\gamma\in \hat K^\times$ satisfies $[\gamma]\in
    A$.  Then, since $\epsilon$ and $\sigma$ commute,
    \begin{equation*}
        [N_{\hat K/\hat F}(\epsilon(\gamma))]_{\hat F} =
        [\epsilon(N_{\hat K/\hat F}(\gamma))]_{\hat F} =
        \epsilon([N_{\hat K/\hat F}(\gamma)]_{\hat F}) = [N_{\hat
        K/\hat F}(\gamma)]_{\hat F}^t.
    \end{equation*}
    Hence $e(\epsilon([\gamma]))=t\cdot e([\gamma])$, and we then
    calculate that $e(T[\gamma])=e([\gamma])$. Therefore
    $e((I-T)[\gamma])=0$, as desired.

    Now since $\Fp[G]$ is an Artinian principal ideal ring, every
    $\Fp[G]$-module decomposes into a direct sum of cyclic
    $\Fp[G]$-modules \cite[Theorem~6.7]{SV}.  Since cyclic
    $\Fp[G]$-modules are indecomposable, we have a decomposition of
    $J=J_\epsilon\oplus J_\nu$ as a direct sum of indecomposables.
    From Section~\ref{se:gm} we know that each of these
    indecomposable modules are self-dual and local, and therefore
    they have local endomorphism rings.  By the
    Krull-Schmidt-Azumaya Theorem (see \cite[Theorem 12.6]{AF}), all
    decompositions of $J$ into indecomposables are equivalent. (In
    our special case one can check this fact directly.)

    On the other hand, we know by \cite{MSS} several properties of
    $J$, including its decomposition as a direct sum of
    indecomposable $\Fp[G]$-modules, as follows. By
    \cite[Theorem 2]{MSS},
    \begin{equation*}
        J=X\oplus \bigoplus_{i=0}^n Y_i,
    \end{equation*}
    where each $Y_i$ is a direct sum, possibly zero, of
    $\Fp[G]$-modules isomorphic to $M_{p^i}$, and $X=N_{\chi}$ for
    some $\chi\in \hat K^\times$ such that $N_{\hat K/\hat
    F}(\chi)\in a^w\hat F^{\times p}$ for some $w$ relatively prime
    to $p$. Moreover, $X\simeq M_{p^r+1}$ for some $r\in \{-\infty,
    0,\dots,n-1\}$.  We deduce that $e(\chi)\neq 0$ and that $e$ is
    surjective on $X$.  Furthermore, considering each $Y_i$ as a
    direct sum of indecomposable modules $M_{p^i}$, we have a
    decomposition of $J$ into a direct sum of indecomposable
    modules.

    We deduce that every indecomposable $\Fp[G]$-submodule appearing
    as a direct summand in $J_\epsilon$ is isomorphic to $M_{p^i}$
    for some $i\in \{0,\dots,n\}$, except possibly for one summand
    isomorphic to $M_{p^r+1}$.  Moreover, we find that $e$ is
    nontrivial on $J_\epsilon$, as follows.  From the hypothesis
    that either $p>2$ or $n>1$ we deduce that $p^r+1<p^n$. Therefore
    since $N_\chi\simeq M_{p^r+1}$ we have $[\chi]\in A$.  Let
    $\theta, \omega\in \hat K^\times$ satisfy $[\theta]=T[\chi]\in
    J_\epsilon$ and $[\omega]=(I-T)[\chi]$.  From $e((I-T)A)=\{0\}$
    we obtain $e(\omega)=0$.  Therefore $e(\theta)\neq 0$.  Observe
    that $\rho^{p^r+1}[\theta]=[1]$.

    We next claim that $e$ is trivial on any $\Fp[G]$-submodule $M$
    of $J_\epsilon$ such that $M\simeq M_j$ for $j<p^r+1$.  Suppose
    not: $M$ is an $\Fp[G]$-submodule of $J_\epsilon$ isomorphic to
    $M_j$ for some $j<p^r+1$ and $e(M)\neq \{0\}$.  Then
    $M=N_\gamma$ for some $\gamma\in \hat K^\times$.  Since $e$ is
    an $\Fp[G]$-homomorphism and $M$ is generated by $[\gamma]$, we
    have $e(\gamma)\neq 0$. But \cite[Proposition 7 and Theorem
    2]{MSS} tells us that $c=p^r+1$ is the minimal value of $c$ such
    that $\rho^{c}[\beta]=[1]$ for $\beta\in \hat K$ with $N_{\hat
    K/\hat F}(\beta)\not\in \hat F^{\times p}$.  Hence we have a
    contradiction.

    Because $J_\epsilon$ decomposes into a direct sum of cyclic
    $\Fp[G]$-modules, we may write $\theta$ as an $\Fp[G]$-linear
    combination of generators of such $\Fp[G]$-modules, and we will
    use this combination and the fact that $e(\theta)\neq 0$ to
    prove that there exists a summand isomorphic to $M_{p^r+1}$ on
    which $e$ is nontrivial.  Let $M=N_\delta$ be an arbitrary
    summand of $J_\epsilon$.  Then $M\simeq M_{j}$ for some $j$. Let
    $[\theta_\delta]$ be the projection of $[\theta]$ on $M$. Since
    $\rho^{p^r+1}[\theta]=[1]$, we deduce that
    $\rho^{p^r+1}[\theta_\delta] = [1]$.  Now if $j>p^{r}+1$ then
    $[\theta_\delta]$ lies in a proper submodule of $M$.  Because
    $\rho M$ is the unique maximal ideal of $M$ and $e$ is an
    $\Fp[G]$-module homomorphism, $e(\theta_\delta)=0$.  On the
    other hand, if $j<p^r+1$ then we have already observed that
    $e(M)=\{0\}$. From $e(\theta)\neq 0$ we deduce that there must
    exist a summand isomorphic to $M_{p^r+1}$ and on which $e$ is
    nontrivial.  Let $U$ denote such a summand.

    Now let $\{V_{\alpha}\}$, $\alpha\in \Ac$, be the collection of
    summands of $J_\epsilon$ apart from $U$.  Hence $J_\epsilon = U
    \oplus_{\alpha\in \Ac} V_\alpha$. Since every summand of
    $J_\epsilon$ is isomorphic to $M_{p^i}$ where $i\in
    \{0,1,\dots,n\}$, except possibly for one summand isomorphic to
    $M_{p^r+1}$, we have (\ref{it:j1}).  From the last paragraph, we
    have (\ref{it:j2}) and (\ref{it:j3}).  Finally, since $e$ is
    trivial on $\Fp[G]$-submodules isomorphic to $M_j$ with
    $j<p^r+1$, we have (\ref{it:j4}).
\end{proof}

\section{Proof of Theorem~\ref{th:main}}\label{se:proof}
\begin{proof}
    We first consider the case $\chr F\neq p$.

    Suppose that $L/F$ is a Galois extension with group $M_{p^i+c}
    \rtimes G$, where $0\le i<n$ and $1 \le c < p^{i+1}-p^i$.  Let
    $K=\Fix_L (M_{p^i+c})$ and identify $G$ with $\Gal(K/F)$.
    Define $\hat F$, $\hat K$, $J$, $J_\epsilon$, and $A$ as in
    Sections~\ref{se:kummer} through \ref{se:jepsilon}. By the
    Kummer correspondence of Section~\ref{se:kummer} and
    Proposition~\ref{pr:kummer}, $L=L_M$ for some $\Fp[G]$-submodule
    $M$ of $J_\epsilon$ such that $M\simeq \Gal(L/K)\simeq
    M_{p^i+c}$ as $\Fp[G]$-modules.  Let $\gamma\in \hat K^\times$
    be such that $M=N_\gamma$.  Since $p^i+c<p^n$, we see that
    $M\subset A$ and so $e$ is defined on $M$. By
    Proposition~\ref{pr:groupoflm} and the discussion following it,
    from $\Gal(L/F)\simeq M_{p^i+c} \rtimes G$ we deduce
    $e(\gamma)=0$.

    Observe that if $p=2$ then from $p^i + c < p^{i+1}$ and $1 \le
    c$ we see that $i > 0$ and hence $n > 1$.  By
    Theorem~\ref{th:jepsilon}, $J_\epsilon$ has a decomposition into
    indecomposable $\Fp[G]$-modules
    \begin{equation*}
        J_\epsilon = U\oplus \bigoplus_{\alpha\in \Ac} V_\alpha
    \end{equation*}
    such that each indecomposable $V_\alpha$ is isomorphic to
    $M_{p^j}$ for some $j\in \{0,\dots,n\}$, $U\simeq M_{p^r+1}$ for
    some $r\in \{-\infty,0,\dots,n-1\}$, $e(U)=\Fp$, and
    $e(V_\alpha)= \{0\}$ for all $V_\alpha\simeq M_{p^i}$ with $0\le
    i\le r$.  Let $U=N_\chi$ for some $\chi\in \hat K^\times$.  Then
    $e(\chi)\neq 0$.

    Because $\rho^{p^i+c-1} M\neq \{0\}$ we know that $J_\epsilon$
    is not annihilated by $\rho^{p^i+c-1}$.  Therefore either
    $\rho^{p^i+c-1}$ does not annihilate $U\simeq M_{p^r+1}$, whence
    $p^r+1\ge p^i+c$, or $p^r+1<p^i+c$ and there exists an
    indecomposable summand isomorphic to $M_{p^j}$ for some $j>i$.

    Suppose first that $p^r+1<p^i+c$ and $J_\epsilon$ contains an
    indecomposable summand $V$ isomorphic to $M_{p^j}$ for some
    $j>i$.  If $j=n$ then by Proposition~\ref{pr:kummer} there
    exists a Galois extension $L_V/F$ such that $\Gal(L_V/K)\simeq
    M_{p^n}\simeq \Fp[G]$.  By
    Proposition~\ref{pr:groupoflm}(\ref{it:lm2}), we have
    $\Gal(L_V/F)\simeq \Fp[G]\rtimes G$. Since $M_{p^{i+1}}\rtimes
    G$ is a quotient of $\Fp[G]\rtimes G$, we deduce that
    $M_{p^{i+1}}\rtimes G$ is a Galois group over $F$.

    If instead $j<n$, then let $\gamma\in \hat K^\times$ such that
    $V=N_\gamma$. Because $e$ is surjective on $U$ we may find
    $\beta\in \hat K^\times$ such that $[\beta]\in U$ and $e(\beta)=
    e(\gamma)$.  Now set $\delta := \gamma/\beta$.  Then
    $e(\delta)=0$ and we consider $N_\delta$.  From $p^j > p^i+c >
    p^r+1$ and $\rho^{p^r+1}[\beta]=[1]$ we deduce that
    $\rho^{p^j-1}[\beta] = [1]$. Then $\rho^{p^j}[\delta] = [1]$
    while $\rho^{p^j-1}[\delta]\neq [1]$, so $N_\delta\simeq
    M_{p^j}$. Let $W=N_\delta$. By Propositions~\ref{pr:kummer} and
    \ref{pr:groupoflm} we obtain a Galois field extension with
    $\Gal(L_W/F)\simeq M_{p^j}\rtimes G$. Since $M_{p^{i+1}}\rtimes
    G$ is a quotient of $M_{p^j}\rtimes G$, we deduce that
    $M_{p^{i+1}}\rtimes G$ is a Galois group over $F$.

    Suppose now that for every $j>i$ there does not exist an
    indecomposable summand isomorphic to $M_{p^j}$.  We claim that
    $r>i$. Suppose not. Then from $p^r+1\ge p^i+c$ we obtain $r=i$
    and $c=1$. Moreover, $U$ is the only summand of $J_\epsilon$ not
    annihilated by $\rho^{p^i}$. Let $\theta\in \hat K^\times$ such
    that $[\theta] = \proj_U \gamma$. If $[\theta]\in \rho U$, then
    $\rho^{p^i}[\gamma]=[1]$, whence $\rho^{p^i}M=\{0\}$, a
    contradiction.  Since $[\theta]\in U\setminus \rho U$ and $\rho
    U$ is the unique maximal ideal of $U$, we obtain that
    $U=N_{\theta}$. Since $e(U)=\Fp$, we deduce that $e(\theta)\neq
    0$.  Now if $V_\alpha\simeq M_{p^j}$ for $j\le r$ then
    $e(V_\alpha) = \{0\}$. Hence $e(V_\alpha)=\{0\}$ for all
    $\alpha\in \Ac$. We deduce that $e(\gamma)\neq 0$, a
    contradiction. Therefore $r\ge i+1$.

    Let $\omega=\rho\chi$ and consider $N_{\omega} = \rho N_\chi =
    \rho U$.  We obtain that $e(\omega)=0$ and $N_\omega\simeq
    M_{p^{r}}$.  By Propositions~\ref{pr:kummer} and
    \ref{pr:groupoflm}, we have that $\Gal(L_W/F)\simeq
    M_{p^{r}}\rtimes G$ for some suitable cyclic submodule $W$ of
    $J_\epsilon$. Since $M_{p^{i+1}}\rtimes G$ is a quotient of
    $M_{p^r}\rtimes G$, we deduce that $M_{p^{i+1}}\rtimes G$ is a
    Galois group over $F$.

    Finally we turn to the case $\chr F=p$.  Recall that we denote
    $M_j \rtimes G$, $j = 1,\dots,p^n$, by $H_j$.  We have short
    exact sequences
    \begin{equation*}
        1\to \Fp\simeq \rho^{p^i+c+k}M_{p^i+c+k+1}\rtimes 1 \to
        H_{p^i+c+k+1} \to H_{p^i+c+k}\to 1
    \end{equation*}
    for all $1\le i<n,$ $1\le c<p^{i+1}-p^i$, and $0\le
    k<p^{i+1}-p^i-c$. For all of these, the kernels are central, and
    the groups $H_{p^i+c+k+1}$ and $H_{p^i+c+k}$ have the same rank,
    so the sequences are nonsplit.  By Witt's Theorem, all central
    nonsplit Galois embedding problems with kernel $\Fp$ are
    solvable. (See \cite[Appendix~A]{JLY}.)  Hence if $H_{p^i+c}$ is
    a Galois group over $F$, one may successively solve a chain of
    suitable central nonsplit embedding problems with kernel $\Fp$
    to obtain $H_{p^{i+1}}$ as a Galois group over $F$.
\end{proof}

\section{Acknowledgements}

Andrew Schultz would like to thank Ravi Vakil for his encouragement
and direction in this and all other projects.  John Swallow would
like to thank Universit\'e Bordeaux I for its hospitality during
2005--2006.

\end{document}